\documentclass[twocolumn,10pt]{article}
\usepackage[top=.5in, bottom=.75in, left=.5in, right=.5in]{geometry}
\pdfoutput=1
\usepackage{amsmath,amsfonts,amscd,amssymb}
\usepackage{graphicx}
\usepackage{epstopdf}
\usepackage{overpic}
\usepackage{cancel}
\usepackage{rotating}
\usepackage{url}
\usepackage{caption}
\usepackage{color}
\usepackage{rotating}
\usepackage{multirow}
\usepackage{wrapfig}
\usepackage{mathtools}
\usepackage{subeqnarray}
\usepackage{setspace}
\usepackage{palatino} 
\usepackage{booktabs}

\usepackage{nomencl}  
\newcommand{\sunderb}[2]{
  \mathclap{\underbrace{\makebox[#1]{$ $}}_{#2}}
}
\immediate\write18{makeindex \jobname.nlo -s nomencl.ist -o \jobname.nls}

\usepackage[bottom,flushmargin,hang,multiple]{footmisc}
\usepackage{lipsum}
\newcommand\blfootnote[1]{%
  \begingroup
  \renewcommand\thefootnote{}\footnote{#1}%
  \addtocounter{footnote}{-1}%
  \endgroup
}

\DeclareMathOperator*{\argmax}{arg\rm{}max}
\DeclareMathOperator*{\argmin}{arg\rm{}min}

\title{\huge{Neural-inspired sensors enable sparse, efficient classification of spatiotemporal data}}
\author{\normalsize{Thomas L. Mohren $^{1,2}$, Thomas L. Daniel $^2$, Steven L. Brunton $^1$, and Bingni W. Brunton$^{2,*}$}\\
\footnotesize{$^1$ Department of Mechanical Engineering, University of Washington, Seattle, WA 98195, United States}\\
\footnotesize{$^2$ Department of Biology, University of Washington, Seattle, WA, 98195 , United States}
}
\date{}

\begin{document}

\maketitle
\makenomenclature
\blfootnote{$^*$ Corresponding author (bbrunton@uw.edu).\\ \noindent \textbf{Code:} \url{github.com/tlmohren/Mohren_WingSparseSensors} }


\vspace{-.3in}
\begin{abstract}
Sparse sensor placement is a central challenge in the efficient characterization of complex systems when the cost of acquiring and processing data is high.  
Leading sparse sensing methods typically exploit either spatial or temporal correlations, but rarely both. 
This work introduces a new sparse sensor optimization that is designed to leverage the rich spatiotemporal coherence exhibited by many systems. 
Our approach is inspired by the remarkable performance of flying insects, which use a few embedded strain-sensitive neurons to achieve rapid and robust flight control despite large gust disturbances.  
Specifically, we draw on nature to identify targeted neural-inspired sensors on a flapping wing to detect body rotation. 
This task is particularly challenging as the rotational twisting mode is three orders-of-magnitude smaller than the flapping modes. 
We show that nonlinear filtering in time, built to mimic strain-sensitive neurons, is essential to detect rotation, whereas instantaneous measurements fail. 
Optimized sparse sensor placement results in efficient classification with approximately ten sensors, achieving the same accuracy and noise robustness as full measurements consisting of hundreds of sensors.  
Sparse sensing with neural inspired encoding establishes a new paradigm in hyper-efficient, embodied sensing of spatiotemporal data and sheds light on principles of biological sensing for agile flight control. \\



\end{abstract}



\section{Introduction}
In both living systems and modern technology, there is a tension between gathering vast and increasing quantities of heterogeneous data (e.g., the internet-of-things), versus acquiring targeted data gathered by specialized sensors~\cite{boyd_six_2011, lazer_parable_2014}.
Large numbers of sensors would provide extensive information about the system and its environment but may, in turn, command high energetic costs.
Indeed, big data demands synthesis and significant processing, often to identify which few features of the data are meaningful, particularly when the crucial information is obscured by large, non-relevant signals or noise.
In contrast, each specialized sensor can extract features tailored to the signal, but unanticipated features in the data may be lost.
The tradeoff between flexibility and efficiency relies in part on the relative difficulty of acquiring, transforming, and performing complex computations on the data.
In addition, local computations alleviate expensive data transfers and may reduce the latency of a decision.
Here, we focus on understanding and designing systems with \emph{sparse and efficient} sensing strategies that leverage both correlations in space and dynamics in time.

Recent advances in sparse sensing rely on the observation that many signals in nature exhibit relatively simple, low-dimensional patterns, so that signal reconstruction or classification can be achieved with a small subset of all possible sensors.
In particular, compressed sensing theory states that if the information of a signal $\mathbf{x}$ is sparse in a transformed basis $\mathbf{\Psi}$, then the signal may be reconstructed from relatively few incoherent measurements~\cite{Candes:2006a,Donoho:2006,Baraniuk:2007,Romberg:2008}.
The number of measurements may be further reduced by taking two additional perspectives.
First, if we do not use a universal transform basis (e.g. Fourier, wavelets, etc.) but instead learn $\mathbf{\Psi}$ from training data, sensor selection may be tailored to a specific task~\cite{manohar2017data}.
Second, when only classification is required, reconstruction can be circumvented and the number of measurements needed are orders-of-magnitude fewer still~\cite{proctor2014exploiting}.
Here we use the sparse sensor placement optimization for classification (SSPOC,~\cite{brunton2016sparse}) approach to identify the locations of a few, key strain sensors tailored to inform body rotation.


We turn to flight control in insects as inspiration of a sensing strategy by which temporal and spatial information are combined.
Flying insects are remarkably adept at making rapid and robust corrections to stabilize their body orientation in response to gusts.
This robust flight control relies on multimodal integration of visual and mechanical information; vision is crucial for flight---indeed, insects rarely fly without it---yet the slow timescale of visual processing cannot support the rapid maneuvers observed in free flight~\cite{collett1975visual, theobald2009wide, sponberg2015luminance}.
Insects accomplish this task using mere tens to hundreds of neurons acting as strain sensors located on their bodies~\cite{nalbach1994halteres, sherman2003comparison, sane2007antennal, taylor2007sensory}, despite the complexity of the surrounding fluid dynamics~\cite{bomphrey2017smart}.
Efficient, distributed sensing and computing has also been explored in nature-inspired engineering~\cite{brunton2015closed}; some examples include insights gained from flying insects~\cite{Combes:2001, Faruque2010jtb1, Faruque2010jtb2}, birds and bats~\cite{song2008aeromechanics, Clark2006jfm}, and fish~\cite{Triantafyllou1995sa, Allen2001jfs, Whittlesey2010bb, Leftwich2012jeb}.
In particular, flying insects sense mechanical deflections using neurons associated with mechanosensory structures known as campaniform sensilla on their wings~\cite{dombrowski1991untersuchungen, dickerson2014control} or their halteres~\cite{dickinson1999haltere, huston2009nonlinear, fox2010encoding}, which are structures derived from wings that function as gyroscopes. 
These mechanosensory neurons are implicated in mediating flight posture control~\cite{dickerson2014control} and encode mechanical stimulus features~\cite{pratt2017neural}. 
Even so, they do not resemble typical engineered sensors, as they do not directly report physical measurement quantities.
Instead, mechanosensitive neurons encode physical strain by a transformation that may be summarized as a temporal filter followed by a nonlinear activation function~\cite{fox2008neural, fox2010encoding}; this encoding has been well characterized in animal experiments~\cite{pratt2017neural}.

\begin{figure}[t!]
 \centering
\includegraphics{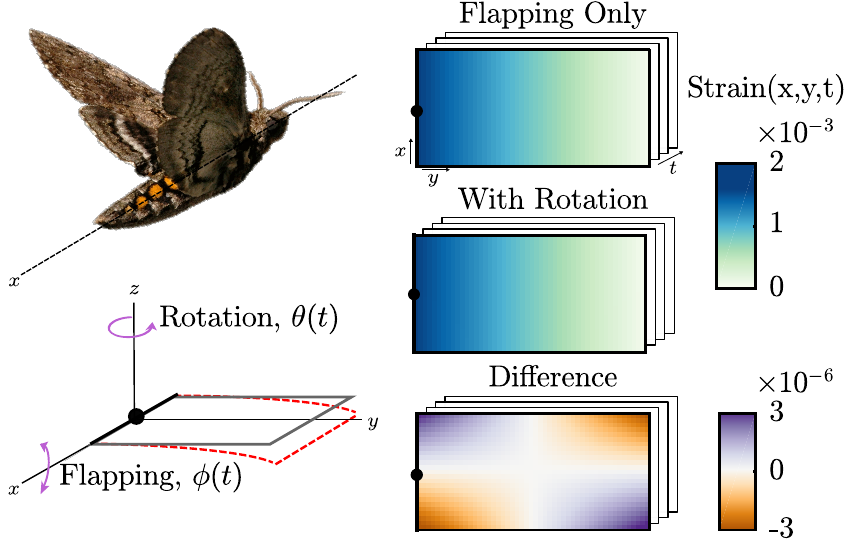}
\caption{
A simulated flapping wing model with and without rotation differ by a twisting mode \emph{three orders-of-magnitude} smaller in magnitude than the dominant flapping mode.
From the flapping wing simulation, we obtain span-wise normal strain over a dense grid on the wing as a function of space $(x, y)$ and time $t$.
Photo of hawkmoth by A. Hinterwirth.
  }
\label{fig:flapper}
\end{figure}

In this paper, we combine SSPOC with mechanical modeling to show that neural inspired encoding of mechanical strain experienced by a flapping wing is necessary for reliable, efficient classification of spatiotemporal data associated with body rotation.
Using raw strain data, it is impossible to distinguish between flapping with and without body rotation.  
This classification task is challenging in part because the spatiotemporal twisting modes induced by body rotation are three orders-of-magnitude smaller than the flapping modes (Fig.~\ref{fig:flapper},~\cite{eberle2015new}).
In contrast, we show that merely 10 neural inspired sensors placed at key locations can achieve similar classification accuracy as a dense grid of sensors distributed over the entire wings.
We find that this performance is robust to large, noisy disturbances added to the biomechanical wing model.
Further, the experimentally derived nonlinear encoder is not unique; instead, exploration of filter function space reveals a large plateau of similar encoders that perform comparably well at this classification task.
Analyzing the locations of these few, key neural inspired sensors offers mechanistic clues of how biology senses in this hyper-efficient regime.

\nomenclature[Rx]{$x$}{Wing chordwise coordinate}%
\nomenclature[Ry]{$y$}{Wing spanwise coordinate}%
\nomenclature[Rz]{$z$}{Wing out of plane axis}%
\nomenclature[Gv]{$\phi(t)$}{Flapping angle}%
\nomenclature[Gh]{$\theta(t)$}{Rotation angle}%

\section{Neural inspired sparse sensors}

Here we take a reverse-engineering perspective to ask: \emph{What is the fewest number of strain sensitive neurons required to inform body rotation, and where should they be placed?}
Answering these questions requires an integrated approach, combining tools from biomechanical simulations, neurophysiology, and sparse optimization.  
This analysis will demonstrate the need for neural-inspired nonlinear filtering in time and the ability to dramatically reduce the number of required sensors through sparse optimization in space.  
All of the code for modeling and classification is openly available and can be found at \url{github.com/tlmohren/Mohren_WingSparseSensors}, and details of our approach are found in the \textit{Supplemental Information}.

%

\begin{figure*}[tb!]
\centering\includegraphics[width=0.85\textwidth]{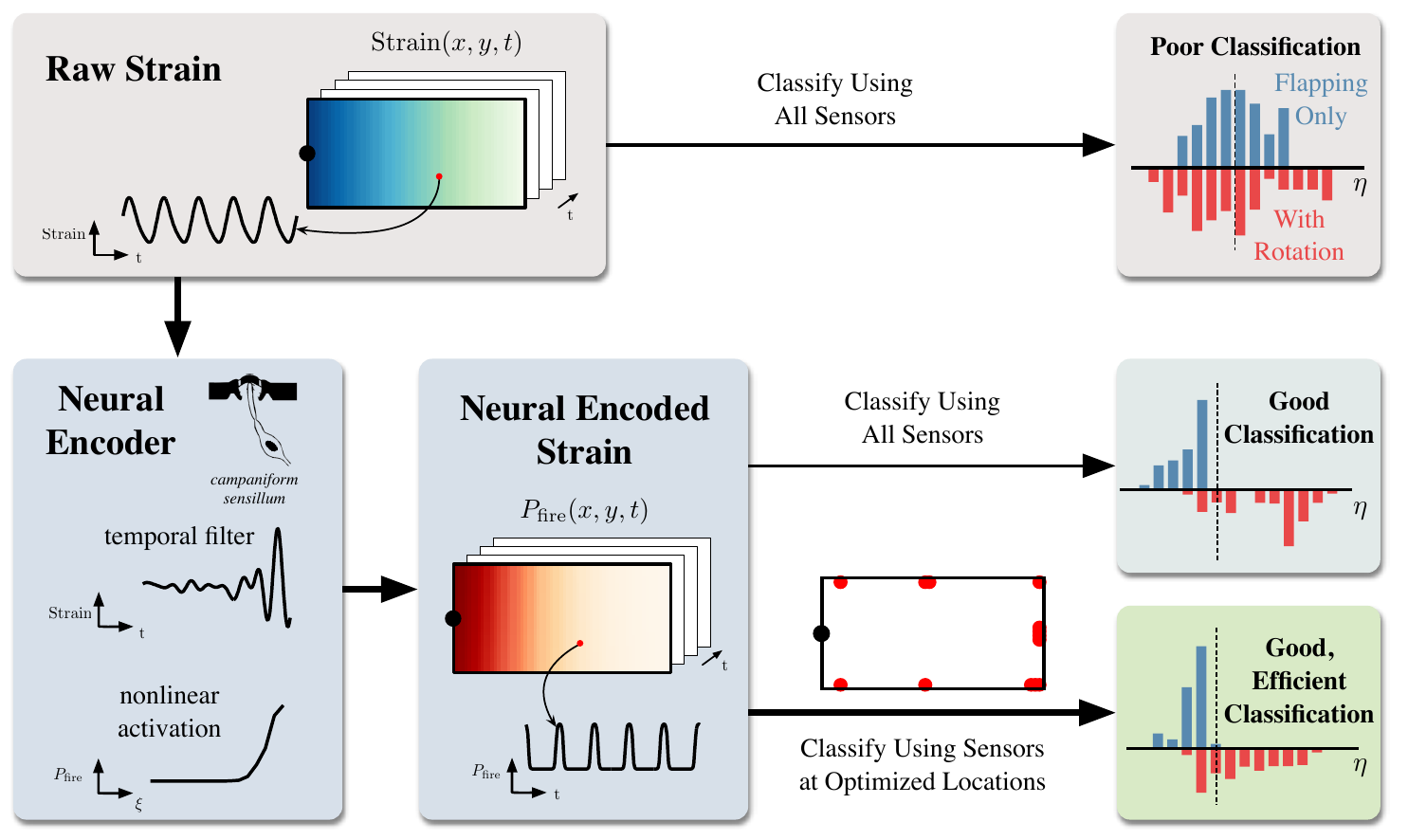}
\caption{ 
A schematic of classifying body rotation using sparse neural inspired strain sensors placed on a flapping wing.
Raw span-wise normal strain is obtained from the structural simulation in two conditions, flapping only and flapping with rotation~\cite{eberle2015new}.
The raw strain in these two conditions are not linearly separable, leading to poor classification even using all the sensor locations.
Alternatively, raw strain is encoded by a neural-inspired filter and transformed into the probability of a mechanosensory neuron firing an action potential~\cite{aljadeff2016analysis, pratt2017neural}.
The neural encoder is approximated by experimental recordings of campaniform sensilla and summarized as a temporal spike triggered average (STA) filter followed by a nonlinear activation function, transforming raw strain into probability of the neuron firing an action potential $P_\text{fire}$.
We define $P_\text{fire}$ to be neural encoded strain.
Neural encoded strain separately well with a linear classifier; further, this performance can be achieved remarkably efficiently using approximately 10 sensors at key locations~ \cite{brunton2016sparse}.
}
\label{fig:neuralEncoding}
\end{figure*}

First, we simulate a flapping wing using an Euler-Lagrange model with parameters based on a hawkmoth~\cite{eberle2015new}.
The flapping wing produces spatiotemporal strain fields sampled at a dense grid on the wing.  
We consider two conditions, given by flapping with and without body rotation (Fig.~\ref{fig:flapper}).
Through a Coriolis force, wing flapping combined with body rotation in an orthogonal axis activates a very small twisting mode in the strain field, and detecting this rotation is a significant challenge. 
We use the simulation data to train a supervised machine learning classifier to distinguish between flapping with and without rotation.  
Random perturbations are added to the flapping and rotational velocities, and the classification accuracy is assessed on validation data from simulations that were not used in training.  

To gauge the role of neural encoding in this task, we compare the performance of classifiers trained using either raw strain from the structural model or neural encoded strain (Fig.~\ref{fig:neuralEncoding}).
The encoding performed by single mechanosensory neurons on the insect wing is approximated by two functions, both of which are derived directly from neurophysiological experiments~\cite{pratt2017neural}.
In short, extracellular recordings of nerve action potentials were made at the wing hinge while mechanical stimuli were delivered to the wing tip through a motor.
Analysis of the mechanical features leading to action potentials were summarized in a temporal spike triggered average (STA) filter followed by a nonlinear activation function~\cite{aljadeff2016analysis}. 
We define the neural encoded strain data as the probability of a mechanosensory neuron firing an action potential.

Next, we solve for the locations of a small subset of sensors among the dense grid on the wing that are sufficient to support classification.
Sensor locations are selected by exploiting the inherent sparsity in the training data.
Our approach uses sparsity-promoting regression and is an extension to SSPOC~\cite{brunton2016sparse}.
Starting with the truncated basis $\mathbf{\Psi}$ and the discriminant vector between the two categories $\mathbf{w}$, we solve for a sparse vector $\mathbf{s} \in \mathbb{R}^n$ that achieves the discrimination $\mathbf{\Psi}^T \mathbf{s} = \mathbf{w}$.
Here, $\mathbf{s}$ has the same shape as the full-state discriminant vector $\mathbf{\Psi} \mathbf{w}$ but contains mostly zeros.
In particular, we use an elastic net penalty to formulate the sparse optimization problem~\cite{tibshirani1996regression, zou2005regularization, clemmensen2011sparse}:
\begin{align}
\mathbf{s} = \argmin_{\mathbf{s'}} \left\lVert \mathbf{w} - \mathbf{\Psi}^T \mathbf{s'} \right\rVert_2 + \alpha \left\lVert \mathbf{s'} \right\rVert_1 + (1 - \alpha) \left\lVert \mathbf{s'} \right\rVert_2,
\label{eq:lossfunction}
\end{align}
where $\lVert \cdot \rVert_2$ is the $\ell_2$ norm, $\lVert \cdot \rVert_1$ is the $\ell_1$ norm, and  $\alpha$ is a hyperparameter of the optimization.
The few non-zero elements of $\mathbf{s}$ correspond to desired sensor locations; these few sensors are able to closely match the performance of full-state classification.

\begin{figure}[t!]
\centering\includegraphics[width=\columnwidth]{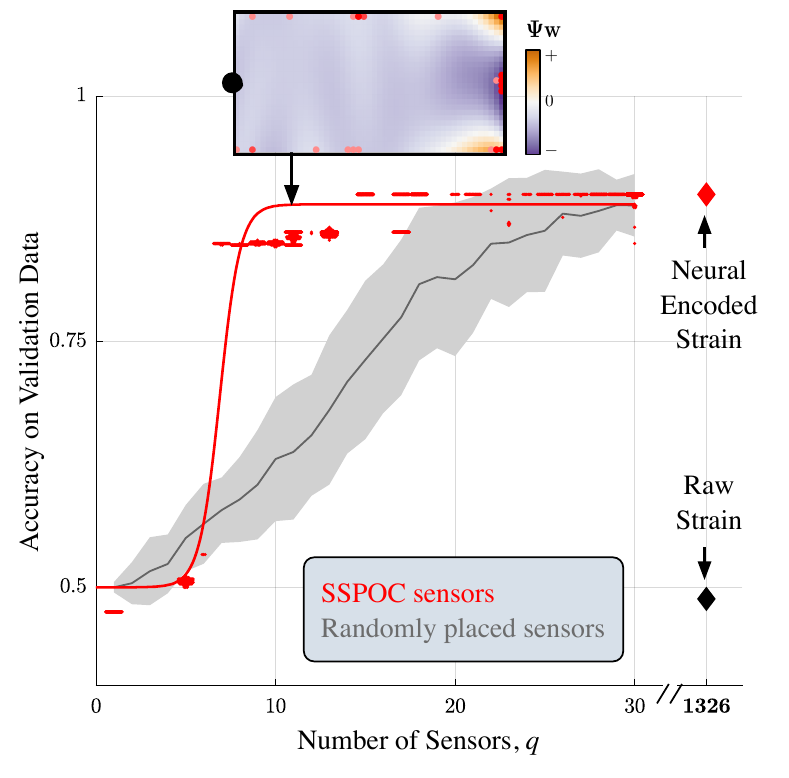}
\caption{
Classification using about 10 neural encoded sensors placed at key locations on the wing achieves accuracy comparable with classification using all sensors.
Flapping wing structural simulations were computed with moderate disturbance amplitudes ($[\dot{\phi}^* ,\dot{\theta}^*] = [0.31,0.1]$ rad/s).
The classification accuracies shown are validated on an epoch of the simulation separate from what had been used for training. 
Sparse sensors are learned from training data from trials with random disturbances using SSPOC (red, each dot is an individual trial) and compared to randomly placed sensors (gray, mean and stdev).
The red line is a sigmoidal fit to the SSPOC sensors accuracy.
The inset shows a probability distribution of SSPOC sensor locations on the wing for $q = 11$ sensors, averaged over 100 training sets with random instances of noisy disturbances.
The opacity of the red dots are proportional to the likeliness of sensor solutions at that location; most sensors are found at the periphery of the wing.
The background of the inset shows the full-state discriminant vector $\mathbf{\Psi}\mathbf{w}$.
 }
\label{fig:Figure_R1}
\end{figure}

\section*{Results}

Our primary result is that classification of flapping with and without body rotation requires neural-inspired encoding of strain data. 
In addition, only a few neural-inspired sensors are needed for classification, showing remarkable robustness to large magnitude disturbances. 
We further characterize how well a family of neural-inspired encoders, including the one derived directly from experimental recordings, are able to perform this classification.

\subsection*{Neural inspired encoders are essential}
The raw strain data reveals that body rotation orthogonal to the axis of flapping introduces a torsional mode in the flapping wing orthogonal to the axis of flapping (Fig.~\ref{fig:Figure_twist},~\cite{eberle2015new}).
Although a signature of the rotation is measurable by strain, the magnitude of this torsional mode is three orders of magnitude smaller than the flapping modes (Fig.~\ref{fig:flapper}). 
In addition, raw strain data of flapping only versus flapping with rotation are not linearly separable.
The data under both conditions are overlapping sinusoidal timeseries at the wing flap frequency; therefore, it is not possible to construct any linear hyperplane that separates them.
A classifier trained on raw strain does no better than chance, even using all available sensors on the wing (Fig.~\ref{fig:Figure_R1}, black diamond).

In contrast, neural encoded strain (Fig.~\ref{fig:neuralEncoding}) enables a linear classifier to detect body rotation, achieving accuracy on validation data of 90\% (Fig.~\ref{fig:Figure_R1}, red diamond).
A spike-triggered average (STA) temporal filter selects a short time-history of raw strain that matches the activation of strain-sensitive wing mechanoreceptors, and a nonlinear activation function transforms the raw strain into a probability of firing an action potential, which we define to be the neural encoded strain.

 \begin{figure}[t]              
\centering\includegraphics[width=\columnwidth]{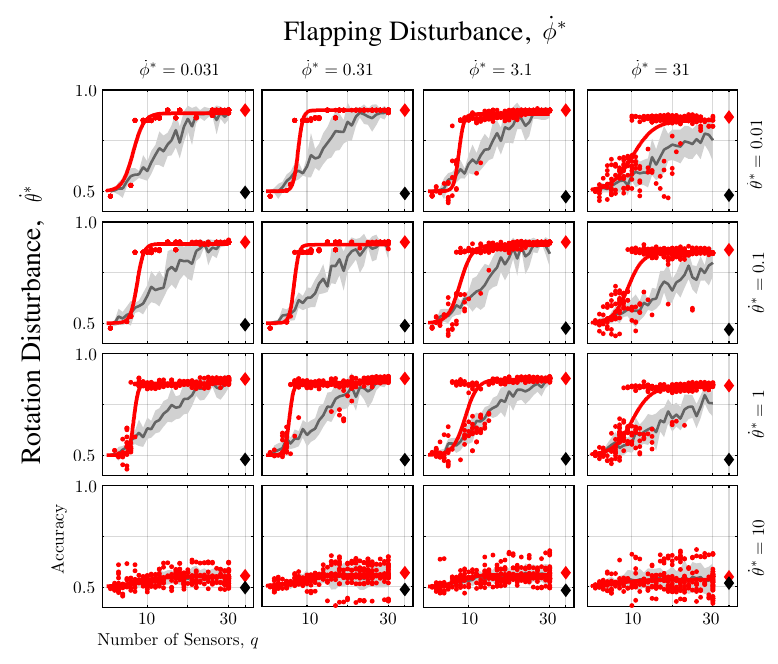}
\caption{
Classification accuracy is robust for moderate to large magnitude disturbances in flapping $\phi$ and in rotation $\theta$.
Each panel of the 4 by 4 grid shows the classification accuracy for varying number of sensors (grey: random sensors, red: SSPOC sensors, black diamond: all sensors without encoding, red diamond: all sensors with encoding). 
The levels of rotation disturbance represent 0.1,1,10, and 100\% of the standard deviation of steady flapping $\bar{ \dot{\phi}}$ and of the magnitude of constant rotation $\bar{\dot{\theta}}$.
}
\label{fig:Figure_R2}
\end{figure}

\subsection*{A few key neural-inspired sensors are required}
\begin{figure}[h]
\centering
\includegraphics{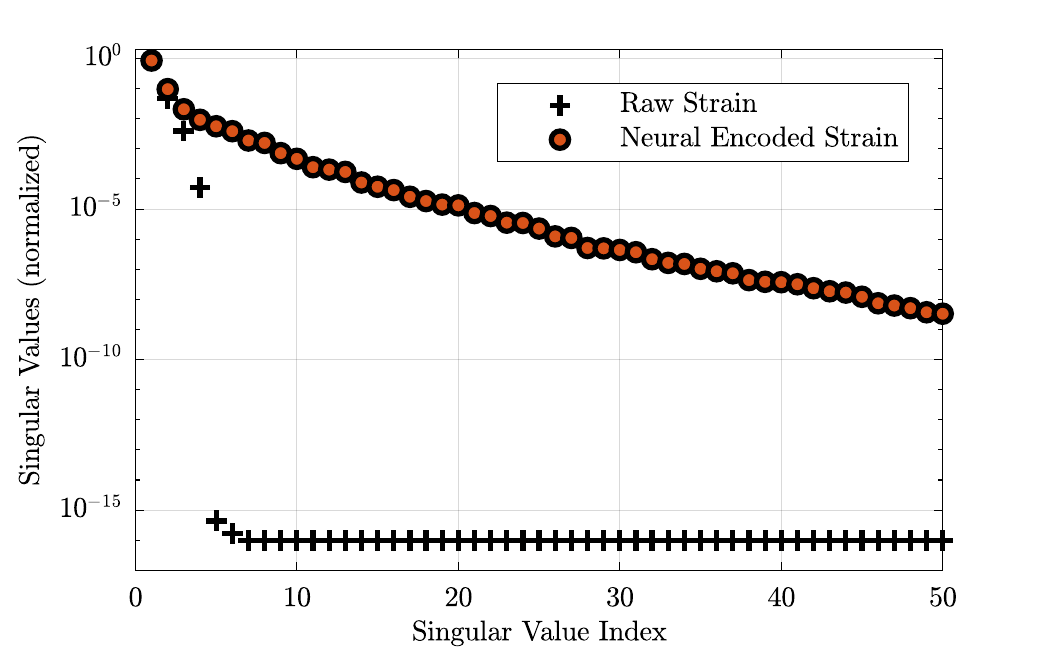}
\caption{The normalized singular values for raw strain (black plus) and neural encoded strain (red circle). }
\label{fig:singular_values}
\end{figure}

\begin{figure}[h]
    \centering
        \includegraphics[width=0.9\linewidth]{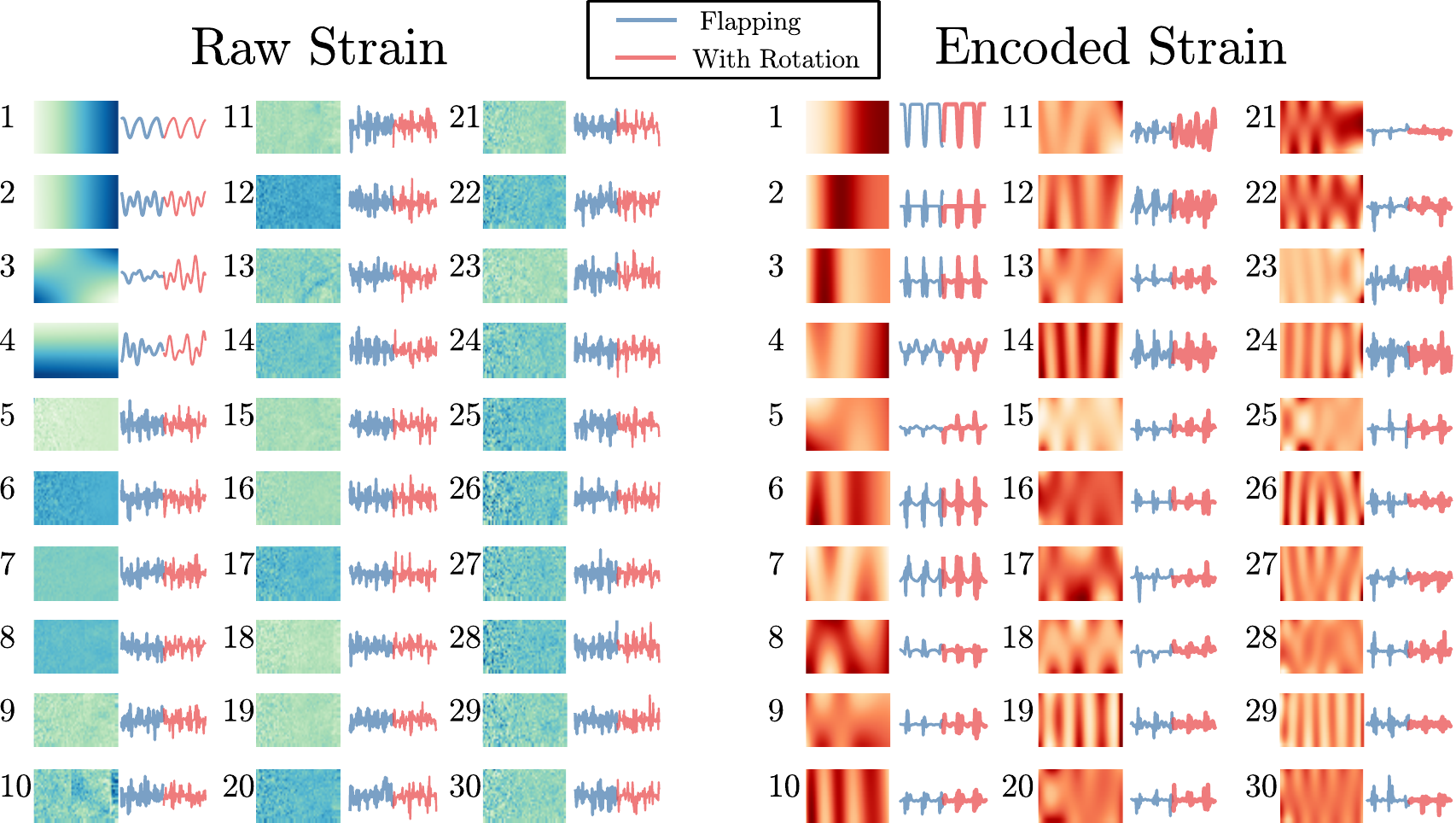}
        \caption{Singular Value Decomposition modes for raw strain(left) and encoded strain (right). The plate shows the mode shape and it's associated number indicates the Singular Value Index. The blue signal shows the presence of that mode when the wing is flapping, the red signal is the presence of that mode for flapping with rotation.}
     \label{fig:svdmodes}
\end{figure}

Importantly, very few of the neural-inspired sensors in the simulation are required for classification, achieving accuracy approaching what is possible with all sensors.
This performance is made possible by exploiting the inherent low rank structure of the data, which is evident in the singular value spectrum of neural encoded strain (Fig.~\ref{fig:singular_values}).
Although the raw strain data is is even lower rank than the neural encoded strain, flapping with and without body rotation remain not linearly separable (Fig.~\ref{fig:svdmodes}).

Indeed, $\sim$25 randomly placed sensors perform just as well on average as using all 1326 sensors (Fig.~\ref{fig:Figure_R1}, grey curve).
It is possible to further reduce the number of sensors by selecting optimized locations, and $\sim$10 SSPOC sensors achieve comparable performance (Fig.~\ref{fig:Figure_R1}, red dots).
The relationship between the number of SSPOC sensors $q$ and validated accuracy follows a sigmoidal shape (Fig.~\ref{fig:Figure_R1}, red curve).

The optimized sensor locations are shown as an inset in Fig.~\ref{fig:Figure_R1} for $q=11$ sensors, and they are distributed at distinct locations at the periphery of the wing.
These locations include the far edge of the wing away from the body, where the full-state discriminant vector $\mathbf{\Psi}\mathbf{w}$ has large amplitude.

 \begin{figure*}[t!]
\centering\includegraphics[width=0.95\textwidth]{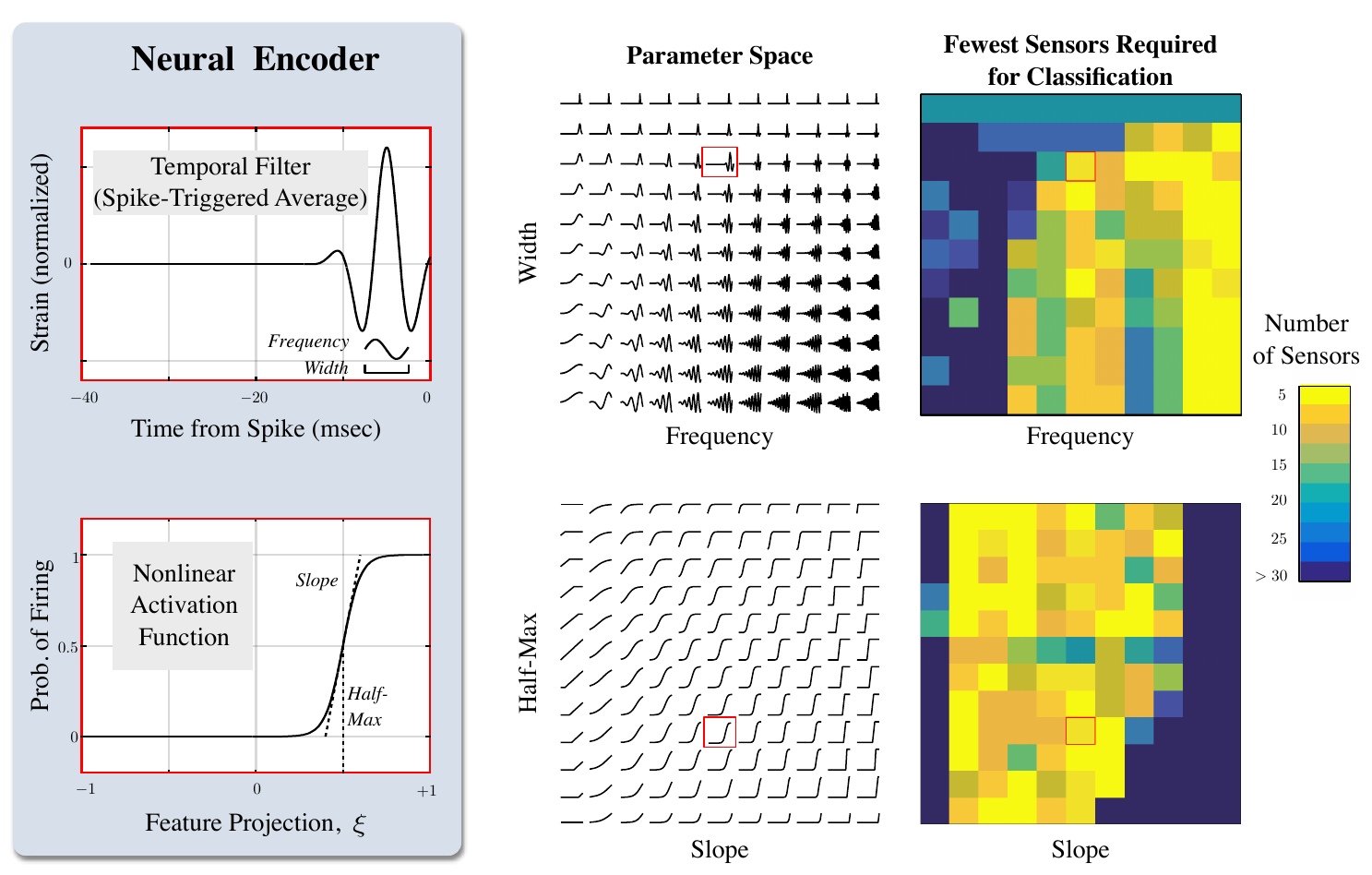}
\caption{
The experimentally derived neural encoders are found on a large plateau in parameter space with similar neural-inspired encoders.
The temporal filter STA is parameterized by its frequency and width (Eq.~\ref{eq:STA}), and the nonlinear activation function is a sigmoid parameterized by its slope and half-max (Eq.~\ref{eq:NLD}).
The top and bottom rows show systematic variations of the STA and the nonlinear activation function by manipulating their parameters, respectively.
The middle column visualizes the family of these neural-inspired encoders.
We assess each encoder by the fewest sensors required to achieve 75\% classification accuracy, and these numbers of sensors are shown as heat maps in the right column.
The experimentally fit encoder functions (red boxes) are well suited to achieve classification along with a family of similar neural-inspired encoders.
  }
\label{fig:Figure_R3}
\end{figure*}

\subsection*{Classification is robust to disturbances}

\begin{figure}[h]
\centering
\includegraphics{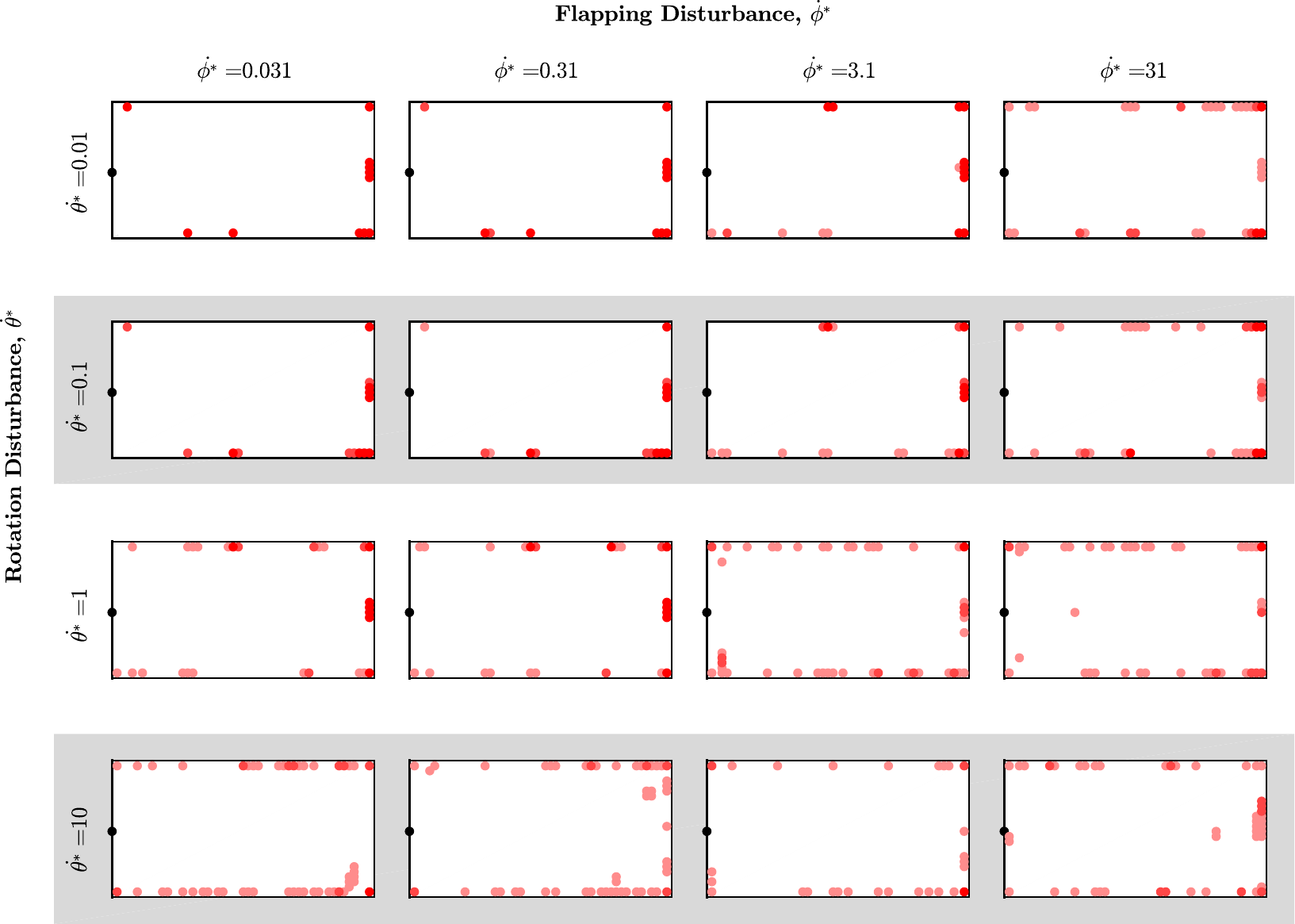}
\caption{The sensor locations for $q=11$ under different disturbance level combinations, darker red indicating a higher probability of sensor placement. }
\label{fig:sensorloc_positive_negative}
\end{figure}

The few key sensors discovered by the SSPOC optimization reliably classify body rotation even when the magnitude of disturbances are large.
Fig.~\ref{fig:Figure_R2} shows the validated classifier accuracies for increasing disturbances in both the flapping $\phi$ and rotational $\theta$ axes.
Smaller disturbances support classification with fewer sensors. 
Even so, the asymptotic full-state accuracy is approached for all rotational disturbances less than 10 radians/sec, at which the disturbances equal the steady rotation velocity.
The performance of sparse sensors is characterized for finer resolutions of disturbances in Fig.~\ref{fig:disturbance_wo_STA}, and the probability distribution of sensors for each disturbance level is shown in Fig.~\ref{fig:sensorloc_positive_negative}.

\begin{figure}[h]
\centering
\includegraphics{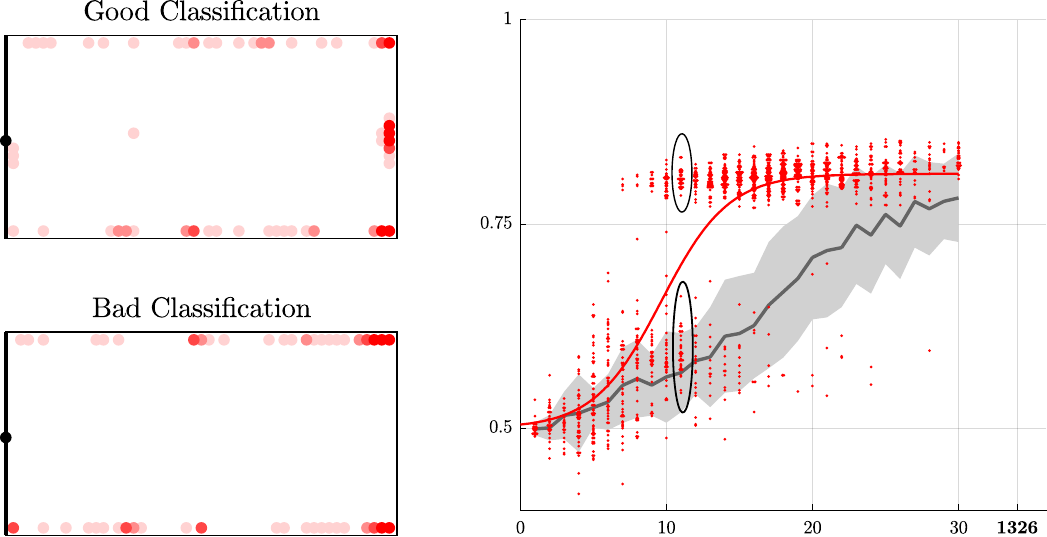}
\caption{The sensor placement for $q = 11$ under disturbance conditions $\dot{\phi}^* = 31$ $\dot{\theta}^* =1$. 
The right figure shows the bimodal distribution of the classification accuracy. 
The sensor locations are shown for the succesful classification group (top) and the poor classification group (bottom), with darker red indicating a higher probability of sensor placement.  }
\label{fig:sensorloc_good_bad}
\end{figure}

Interestingly, when the mean classification degrades for larger disturbances, the distribution of accuracy at a given number of sensors $q$ becomes bimodal. 
In other words, sparse sensor optimization on some sets of training data achieve accuracy that approaches the asymptotic full-state accuracy, but other random instances lead to poor classification.
Comparing the sensor location distributions for the good classification versus the poor classification cases, we see that sensors at the far edge of the wing away from the body are crucial for classification (Fig.~\ref{fig:sensorloc_good_bad}).

\subsection*{Variations on the theme of experimentally derived neural encoders}

So far, we have used a parameterized neural encoder fit directly to electrophysiological recordings of campaniform sensilla in insects~\cite{pratt2017neural}. 
Now, we explore the effects of systematic variations to the neural encoder's parameters to determine whether the experimentally derived encoder is uniquely suited to the task.
The temporal filter and the nonlinear activation function both have two parameters each.
We vary each pair of parameters while holding the others fixed at their experimentally derived values.

The performance achievable by this family of neural encoders is summarized by the fewest sensors required to achieve 75\% classification accuracy.
For each encoder, a full sweep of validated accuracy is computed with at least 10 iterations of random disturbance at each value of $q$.
A sigmoidal fit of the relationship between $q$ and accuracy (as in the red curve in Fig.~\ref{fig:Figure_R1}) is then used to determine at what $q$ the accuracy exceeds 75\%.
For some regimes in the encoder parameter space, this accuracy is never achieved for any number of sensors.

\begin{figure}[h]
 \centering\includegraphics[width=\linewidth]{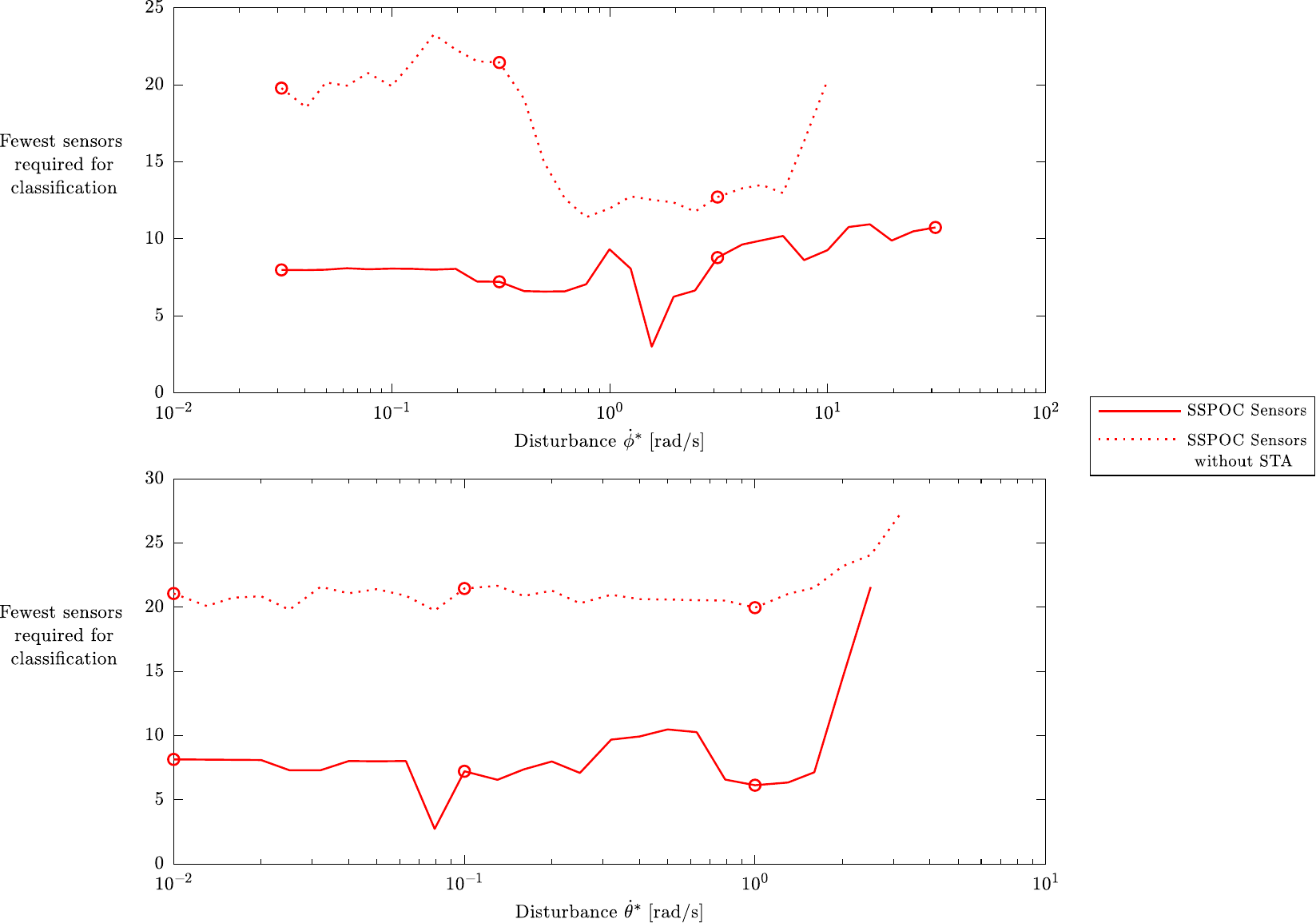}
\caption{
The top plot shows a more detailed plot of the number of sensors, q, required for 75\% accuracy versus an increase in flapping disturbance $\dot{\phi}^*$, with $\dot{\theta}^*$ constant at 0.1. 
The dotted line shows classification accuracy without STA. 
The circle represent the number of sensors for the plots that were shown in figure \ref{fig:Figure_R2}, corresponding to the red bar in the matrix figure on the top right. 
The bottom plot shows a detailed plot for constant $\dot{\phi}^*=0.31$ and varying $\dot{\theta}^*$.
}
\label{fig:disturbance_wo_STA}
\end{figure}

The temporal filter STA has two parameters, frequency and width.
The top row of Fig.~\ref{fig:Figure_R3} shows that the experimentally derived STA (in red boxes) is surrounded by a large plateau in parameter space with comparably STA-like functions.
Further, higher frequency filters tend to perform better, whereas the width of the filter is less crucial as long as it is not too narrow.
In the limit of the narrowest STA, the temporal filter acts as an identity and does not transform the data; in other words, here the encoding is achieved by the nonlinear activation function alone.
The fact that this regime of parameter space is still able to classify rotation, albeit requiring a larger number of sensors, hints at the importance of the nonlinearity.
The STA acts as a temporal filter for disturbances, and without it, classification accuracy degrades for larger noise amplitudes (Fig.~\ref{fig:disturbance_wo_STA}).

Campaniform sensilla nonlinear activation functions generally have a sigmoidal shape (although variations have been observed in experiments~\cite{pratt2017neural, fox2010encoding}).
The bottom row of Fig.~\ref{fig:Figure_R3} shows that the half-max of the sigmoidal function does not impact classification accuracy.
Similarly, the precise slope of the sigmoid is not crucial, as long as it is not too sharp or too shallow.
In the limit of unit slope with zero half-max (middle of left-most column of parameter space), the nonlinear activation function becomes linear.
Without this nonlinearity, classification never achieves 75\% accuracy; in other words, the nonlinear activation function is \emph{required} for classification.

Although the experimentally derived neural encoders are well suited to perform body rotation classification, they are not unique.
In the context of this nature inspired classification task, the observed properties of campaniform sensilla are found in a large parameter space of similar encoders, most of which are able to support robust and sparse classification of body rotation.

\section*{Discussion}

This paper takes inspiration from nature to demonstrate how classification of subtle dynamic regimes in spatiotemporal data can be achieved with remarkably few sensors.
Specifically, we explore how strain sensitive neurons on a flapping wing can detect body rotation, an ethologically relevant task for flying insects.
We show that the task can be accomplished efficiently with very few sensors, even in the presence of large disturbances.
This approach takes advantage of the ability of neurons to encode data with a convolution in time followed by a nonlinear decision function.

The perspectives presented in this paper are related to several prominent domains.
Here we highlight the relationships between neural-inspired sparse sensors and three distinct fields of research, namely deep neural networks, optimal stimulus encoding, and data-driven representation of dynamics.

In the first connection, we note that inspiration from natural neural computation originally gave rise to the study of connectivism and neural networks as an approach in machine learning~\cite{mcculloch1943logical,rumelhart1986learning}.
The recent astonishing success of deep, convolutional neural networks in solving previously intractable problems has relied on the sheer size and complexity of both the networks and the training data~\cite{lecun2015deep,goodfellow2016deep}. 
These deep neural networks have been compared to the abstract, generalized computations performed by the mammalian neocortex~\cite{yamins2016using}.
In contrast, our approach occupies the opposite limit, discovering hyper-efficient solutions to a specific task by learning a minimal set of neural-inspired units.
In addition to neural encoding, our sensors are embedded in a physical simulation, which means they are implicitly leveraging the embodied computation performed by the biomechanical structure itself.

Second, there is a rich body of literature exploring the hypothesis that neural encoding is optimized to efficiently represent input stimulus, usually defined by maximizing mutual information or optimal encoding~\cite{barlow1961possible, sharpee2004analyzing}.
However, here we consider that representation of the stimulus is not an end in itself, but that the animal ultimately gathers information in order to make decisions, act on this information, and control their interactions with the external world.
It follows that the classification framework we have explored here may be embedded in a dynamic, closed-loop control framework, where the sensors inform actuators to interact effectively with a physically realistic environment.

Third, the locations on the wing where sparse neural-inspired sensors are placed (Fig.~\ref{fig:Figure_R1})  do not resemble the locations of campaniform sensilla on a hawkmoth wing~\cite{dombrowski1991untersuchungen,dickerson2014control}.
One difference is that the optimization problem we solve in Eq.~\ref{eq:lossfunction} does not constrain the relative spatial locations of the sensors, whereas an insect's sensors are constrained by biological structures such as the trajectories of the wing veins. 
In addition, despite being challenging, the body rotation detection task we have formulated here is likely only one of the many functions associated with wing mechanosensors.
Interestingly, if we consider that a wing experiences unsteady fluid forces, sampling simultaneously from multiple sensors in space can provide information similar to sampling at a single location in time.
A set of spatial sensors in a fixed configuration can report phase-delayed information, providing natural coordinates suitable for representing complex spatiotemporal dynamical systems~\cite{brunton2016koopman,brunton2017chaos}.

Finally, this work establishes a novel framework for design of hyper-efficient, embodied autonomous sensing.
We envision that the framework motivates development of hardware demonstrations using flexible materials~\cite{majidi2014soft,dean2017robust}.
Recent innovations in 3D printing technology have enabled manufacturing of flexible structures with embedded strain sensors~\cite{muth2014embedded}.
Some of these sensors are capacitive devices with low temporal resolution~\cite{shin2016sensingarticle}, and they have been limited in number and energy budget on small devices.
We suggest that our neural inspired sensing perspective may pivot both of these limitation into advantages in the design of autonomous micro-robotic implementations.

\subsection*{Acknowledgements}
{\small
We are grateful for preliminary analyses by Jared Callaham and Sam Kinn, and for helpful discussion with Bradley Dickerson, Annika Eberle, Nathan Kutz, Krithika Manohar, Eurika Kaiser, and Brandon Pratt.
This work was funded by the AFRL grant (FA8651-16-1-0003) to BWB and SLB;
AFOSR grant (FA9550-18-1-0200) to SLB;
AFOSR grant (FA9550-18-1-0114), the Alfred P. Sloan Foundation, and the Washington Research Foundation to BWB;
AFOSR grant (FA9550-14-1-0398) and the Komen Endowed Chair to TLD.
TLM was supported by the Washington Research Foundation Innovation Graduate Fellowship in Neuroengineering.
}
\subsection*{Author contributions}
{\small
B.W.B. and T.L.D. conceived of the study; T.L.M. carried out the simulations; T.L.M., T.L.D., S.L.B., and B.W.B. analyzed the results; T.L.M., T.L.D., S.L.B., and B.W.B. wrote the paper.
}
%
%

\onecolumn

\clearpage
\newpage
\section*{Supplemental Information}

This section supplements the main text, providing details of our methods as well as additional figures to support our results.
We also include here a nomenclature to summarize the variables and notational conventions we are using in this paper.
All of the code we developed as part of this paper are openly available as a GitHub repository as described in Sec.~\ref{ss:codeanddata}.

First, Sec.~\ref{ss:structuralmodel} describes the flapping wing simulation with and without rotation, its implementation, and extraction of span-wise normal strain from this simulation.
Raw strain computed on a dense grid on the wing is filtered by a neural encoder model.
Our neural encoder model is taken directly from experimental recordings of campaniform sensilla on moth wings, and Sec.~\ref{ss:neuralencoding} briefly describes these experiments as well as how we derived a functional approximation of the neural encoding.
Sec.~\ref{ss:classification} details how raw strain and neural encoded strain data is used to formulate and solve a classification task.
Next, we select a small subset of sensors among the dense grid using sparsity-promoting optimization using an algorithm described in Sec.~\ref{ss:sspoc}, and these few sensors are able to solve the same classification task.
The relationship between number of sensors and validated classifier accuracy is fit by a sigmoidal function (Sec.~\ref{ss:sigmoid}).
Finally, Sec.~\ref{ss:supplresults} include 5 additional figures to supplement the results described in the main text.

\printnomenclature
\newpage

\subsection{Code and Data Access}  \label{ss:codeanddata}

All code we developed to run the simulations and perform the analyses is available in a repository accessible at \url{github.com/tlmohren/Mohren_WingSparseSensors}. 
The code is implemented in MATLAB 2015a. 
To reproduce the figures, a basic MATLAB installation is sufficient. 
To run the simulations and solve the sparse optimization problems, our code is dependent on two toolboxes: MATLAB's symbolic toolbox and CVX (\url{www.cvx.com}).

\subsection{Structural Model} \label{ss:structuralmodel}

To simulate a flapping wing, we use an Euler-Lagrange model for a flapping flat plate and obtain strain for different prescribed inertial rotations.
The model is based on Eberle et al. \cite{eberle2015new} and modified to allow additional velocity disturbances.
The flat plate has a span of 50 mm, chord length of 25 mm, a thickness of 0.0127 mm, and an E-modulus of 3 GPa; these parameters are chosen to be consistent with previous work on hawkmoth structural wing studies \cite{eberle2015new,combes2003flexural2}.

\begin{figure}[t]
 \centering\includegraphics{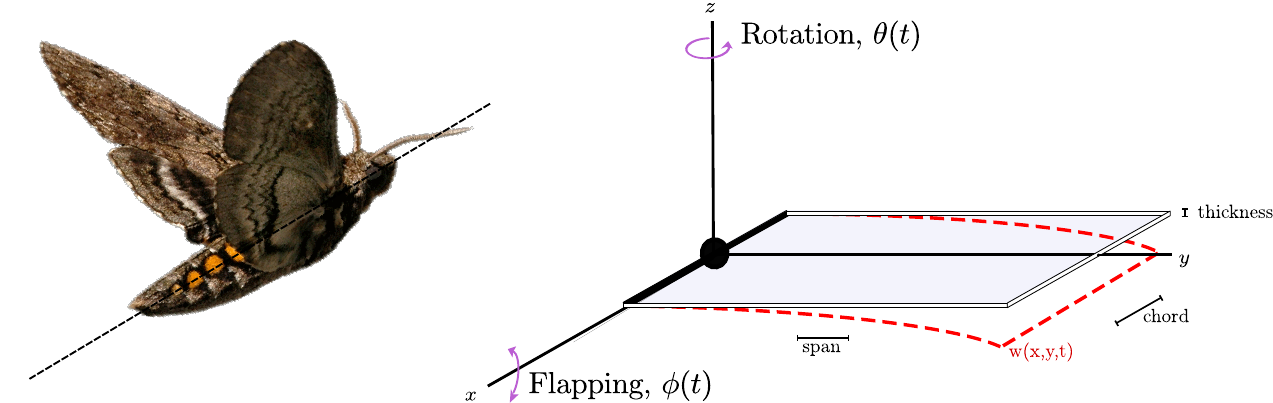}
\caption{
The flapping wing is modeled by a flapping flexing plate, with parameters scaled to match a hawkmoth wing. The plate is excited by the flapping angle $\phi(t)$ around the $x$-axis, and the rotation angle $\theta(t)$ around the $z$-axis.
}
\label{fig:Figure_01B}
\end{figure}

The simulated wing flaps with an amplitude of $\pi/6$ radians at a frequency of $f_\phi = 25$ cycles per second (Hz).
A harmonic at 50 Hz is at $1/5$ the magnitude of the dominant frequency. 
Specifically, the steady flapping is

\begin{align}
\phi(t)    = \frac{\pi}{6} \left( \sin( 2 \cdot 10^{-3} \pi f_\phi t)  + \frac{1}{5}\sin(4  \cdot 10^{-3} \pi f_\phi t) \right),\label{eq:phi}
 \end{align}
 where time $t$ has units of milliseconds.
 
In addition to flapping $\phi(t)$, the wing is perturbed by one of two different inertial rotation velocities $\dot\theta$:
\begin{align}
\dot{\theta} & = 0 \text{ rad/s (without rotation)}, \\
\dot{\theta} & = 10 \text{ rad/s (with rotation)}.
\end{align}

\subsubsection{Physics of the flapping wing} 
A useful perspective on the forces present on the plate can be gained by defining a rotating reference frame $R$. 
We define the position of a point on the centerline of the wing at a distance $L$ from the origin as $\mathbf{r}$ and the rotation of the local frame as $\mathbf{\omega}$:

\begin{align}
\mathbf{r} &= [0,L \cos(\phi), L \sin(\phi) ], \\
\mathbf{\omega} &= [0,0,\dot{\theta} ].
\end{align}

The kinematics in the local frame $R$ is related to global acceleration $I$ as 

\begin{align}
 \left|  \frac{ \partial^2 \mathbf{r}}{\partial t^2}  \right|_I  & =
 \left|    \frac{ \partial^2 \mathbf{r}}{\partial t^2} \right|_R   +     
  \frac{ \partial \mathbf{\omega} }{\partial t}  \times \mathbf{r}  + 
  2   \mathbf{  \omega}\times \frac{ \partial \mathbf{r} }{\partial t}   +  
  \vphantom{  \frac{\partial^2 \bf{r}}{\partial t^2}}  \mathbf{  \omega} \times ( \mathbf{  \omega} \times \mathbf{r} )   \label{eq:coriolis_2},  \\
 \left|  \frac{ \partial^2 \mathbf{r}}{\partial t^2}  \right|_I  & = \begin{bmatrix}
 2 L \sin( \phi ) \dot{\phi} \dot{\theta}   - L\cos (\phi) \ddot{\theta} \\
- L \sin( \phi ) \ddot{\phi }    - L \cos(\phi) ( \dot{\phi}^2 + \dot{\theta}^2 )  \\
 L \cos(\phi) \ddot{\phi}   - L \sin (\phi) \dot{\phi}^2 
 \end{bmatrix}^T    \begin{bmatrix} i \\ j \\ k  \end{bmatrix}. 
\end{align}

\begin{figure}[t]
 \centering\includegraphics[width=0.3\textwidth]{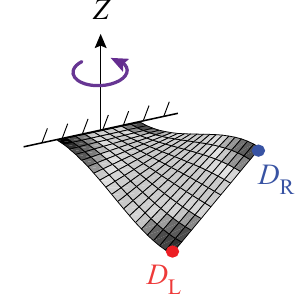}
\caption{A flat plate undergoing both flapping and rotation will undergo a twisting deformation. For the flapping flat plate, this results in an additional strain three orders smaller than strain caused by the main bending deformation mode \cite{eberle2015new}. 
}
\label{fig:Figure_twist}
\end{figure}

Here the accelerations in the $i$ direction are present only for nonzero $\dot{\theta}$. 
Acceleration in this direction results in a twisting mode on the flat plate (Fig.~\ref{fig:Figure_twist}, \cite{eberle2015new}). 
This acceleration is generally referred to as the Coriolis acceleration. 
With use of the small angle approximation, it becomes clear that this acceleration occurs at twice the flapping frequency of $\phi\propto \sin (f_\phi t)$:
\begin{align}
 2 L \sin( \phi ) \dot{\phi} \dot{\theta} &   \approx   2 L \phi  \dot{\phi} \dot{\theta} \\
  & \propto   \sin ( f_\phi t ) \cos(f_\phi t )\\
  & \propto        \sin ( 2 f_\phi t ).
\end{align}

\subsubsection{Adding random disturbances on velocity} 
To simulate structurally relevant noise experienced by flapping wings, we added white noise disturbances to the steady flapping and rotation velocities.
Specifically, the total flapping velocity $\dot\phi_T(t)$ is the sum of steady rotation velocity $\dot\phi(t)$ and the added disturbance $\dot{\phi}^*(t)$.
Similarity, the total inertial rotational velocity $\dot\theta_T(t)$ is the sum of steady rotational velocity $\dot\theta$ and the added disturbance $\dot{\theta}^*(t)$.
 
We modeled disturbances as band limited white noise, summing 15 sinusoids at random frequencies between 1 and 10 Hz and at random phases:

\begin{align}
\dot{\phi}^*(t) & = A_{\dot{\phi}^* }       \sum_{i=1}^{15}    \sin( 2 \rho t + \kappa )  ,   \\
\dot{\theta}^* (t) &  = A_{\dot{\theta}^* }    \sum_{i=1}^{15}  \sin( 2 \rho t + \kappa ) . 
\end{align}

Here the amplitude of the disturbances  are $A_{\dot{\phi}^* }$ and $A_{\dot{\theta}^* }$.
We chose the range of disturbance amplitudes to correspond to 0.1,1,10, and 100\% of the standard deviation of steady flapping $\bar{ \dot{\phi}}$ and 0.1,1,10, and 100\% of the magnitude of constant rotation $\bar{\dot{\theta}}$.
There are two random variables, $\rho$ and $\kappa$;
$\rho$ is a random frequency drawn in the range of $[1, 10]$ Hz, and $\kappa$ is a random phase drawn in the range of $[0, 2\pi]$.

\subsubsection{Ramp up of flapping and rotational velocities in simulation}

\begin{figure}[t]
 \centering\includegraphics{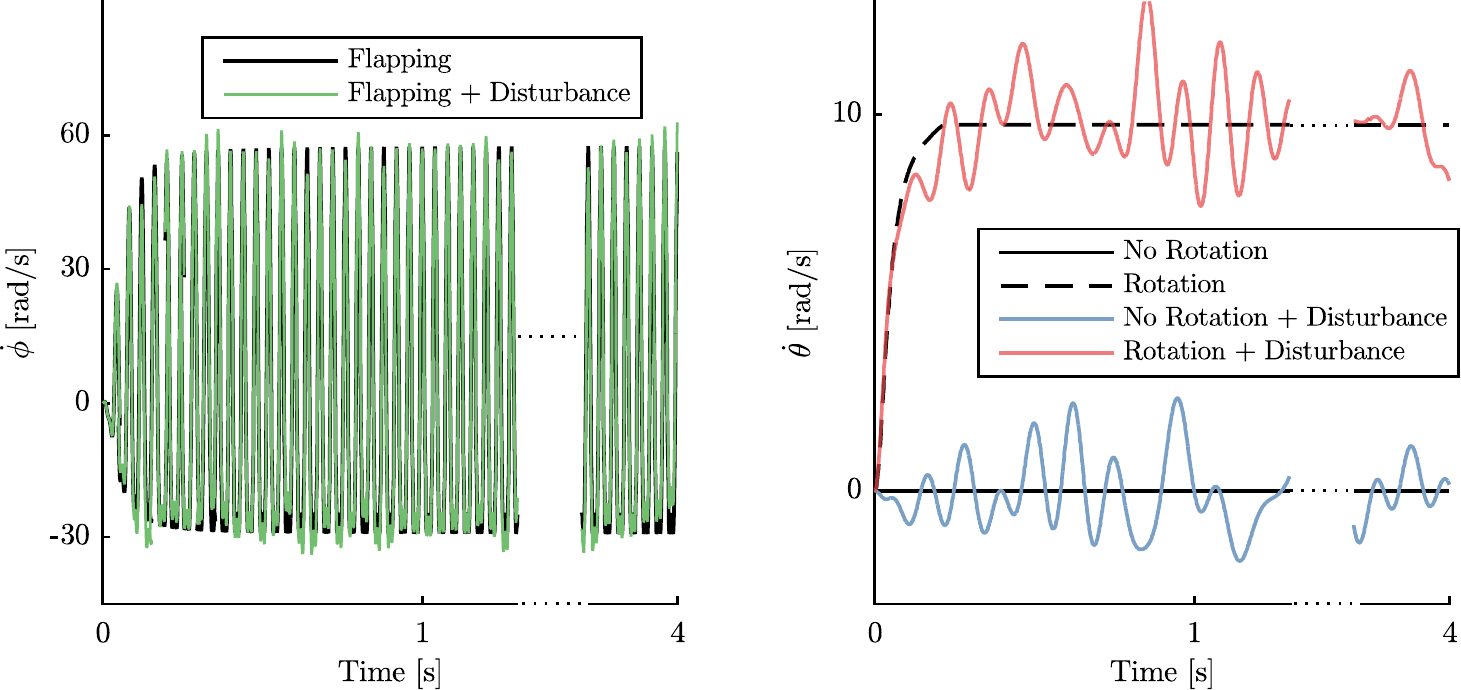}
\caption{
The left graph shows the flapping velocity without disturbance (black) and with disturbance (green, standard deviation $ \sigma = 3.1$ (rad/s)). 
The right graph shows rotational velocity. 
The solid black line shows a rotational velocity of zero, where the blue line shows no steady rotational velocity but a disturbance with $\sigma = 1$ (rad/s). 
The striped black line shows a steady rotational velocity of 10 (rad/s). 
The yellow has a steady rotational velocity of 10 (rad/s) with a disturbance of $\sigma = 1$ (rad/s). 
The simulations had a 960 millisecond transient phase that was discarded for the classification task. 
}
\label{fig:rampup}
\end{figure}

If the Euler-Lagrange simulations were started at the flapping and rotational velocities described above, undesirable high frequency deformation modes would be excited.
To prevent these modes from contaminating the results, both flapping and rotational velocities are multiplied by a sigmoidal startup ramp $\nu$:
\begin{align}
\nu = \frac{ (2\cdot 10^{-3} \pi f t)^3 }{10 +(2\cdot 10^{-3}  \pi f t)^3}. \label{eq:ramping} 
\end{align}

Fig.~\ref{fig:rampup} shows the the total flapping and rotational velocities experienced by our model, including this startup ramp, with and without added disturbances.

\subsubsection{Implementation of model}

We use rectangular shape functions that are fixed in $x$ and $y$ at the corners, but that are free to deform in $w$ at the free edge. 
\begin{align}
\mathbf{r}(x,y,t) = [x,\   y,\   w(x,y,t)]^T 
\end{align}

Using the shape functions, we obtain 3 degrees of freedom for the two free wing corners  $\mathbf{q}$,
\begin{align}
\mathbf{q}(t) = [ \delta_3 ,\  \phi_3,\   \theta_3,\   \delta_4,\   \phi_4,\    \theta_4]^T. \label{eq:q}
\end{align}

The displacement is then a function of the shape functions $\mathbf{N}(x,y)$ and $\mathbf{q}$:
\begin{align}
w(x,y,t) = [ \mathbf{N}(x,y)]^T \mathbf{q}(t).
\end{align}

Using Lagrange's equation, we can then obtain the system of equations:
\begin{align}
\frac{ d^2 \mathbf{q}(t)}{d t^2} = - \mathbf{M}^{-1}\mathbf{M}_a \frac{d v_0}{dt} + \left( \frac{ d \mathbf{\Phi} }{dt} \right)^2 \mathbf{q}(t) - \mathbf{ M}^{-1} \mathbf{K}  \mathbf{q}(t) + \mathbf{M}^{-1}\mathbf{ I}_c \mathbf{\Omega} - \mathbf{M}^{-1} \eta \frac{ d \mathbf{q}(t)}{dt}. 
\end{align}

The system of ODE's from \cite{eberle2015new} with modified rotation angles are solved with MATLAB's ODE45 (5th order Runga-Kutta solver) for $t\in [1,4000]$ milliseconds.

\subsubsection{Computing strain}  

The spanwise strain over the wing $\epsilon_y(x,y,t)$ relates directly to the local curvature through the double partial derivative of the shape function: 

\begin{align}
\epsilon_y(x,y,t) = - \frac{h}{2} \frac{\partial^2 w(x,y,t)}{\partial y^2}. \label{eq:strain_dev}
\end{align}

The chordwise strain $\epsilon_x$ is much smaller than the spanwise strain $\epsilon_y$, and there is no indication that campaniform sensilla can detect shear strain $\epsilon_{xy}$. 
Therefore, we will use the spanwise strain $\epsilon_y$ for our experiments and refer to it as $\epsilon$ in this paper.

Since $w(x,y,t)$ is a continuous function, we can specify our sensor locations anywhere on the wing surface. 
We chose to compute strain over a grid with 0.1 centimeter spacing starting at the edges.
This space results in 51 spanwise and 26 chordwise points, for a total of 1326 sensor locations. 

\nomenclature[Rt]{$t$}{Time}%
\nomenclature[RX]{$ \mathbf{X} $}{Matrix of data gathered at multiple times}%
\nomenclature[Rx]{$ \mathbf{x} $}{Vector of data at one snapshot, measured at sensor locations}%
\nomenclature[RXtr]{$\mathbf{X_\text{train}}$}{Part of the data matrix used to train the classification algorithm}%
\nomenclature[RXte]{$\mathbf{X_\text{test}}$}{Part of the data matrix reserved for testing the classification accuracy}%
\nomenclature[Gv]{$ \dot{\phi} $}{Flapping velocity}%
\nomenclature[Gh]{$ \dot{\theta} $}{Rotational velocity}%

%

\subsection{Neural encoding} \label{ss:neuralencoding}
%

The action potential responses of campaniform sensilla to strain on the hawkmoth wing had been characterized by Pratt et al.~\cite{pratt2017neural}.
Here we briefly describe these experiments and the neural encoding functions fit to the experimental data.

\subsubsection{Experimental electrophysiological recordings}
\begin{figure}[h!]
 \centering\includegraphics{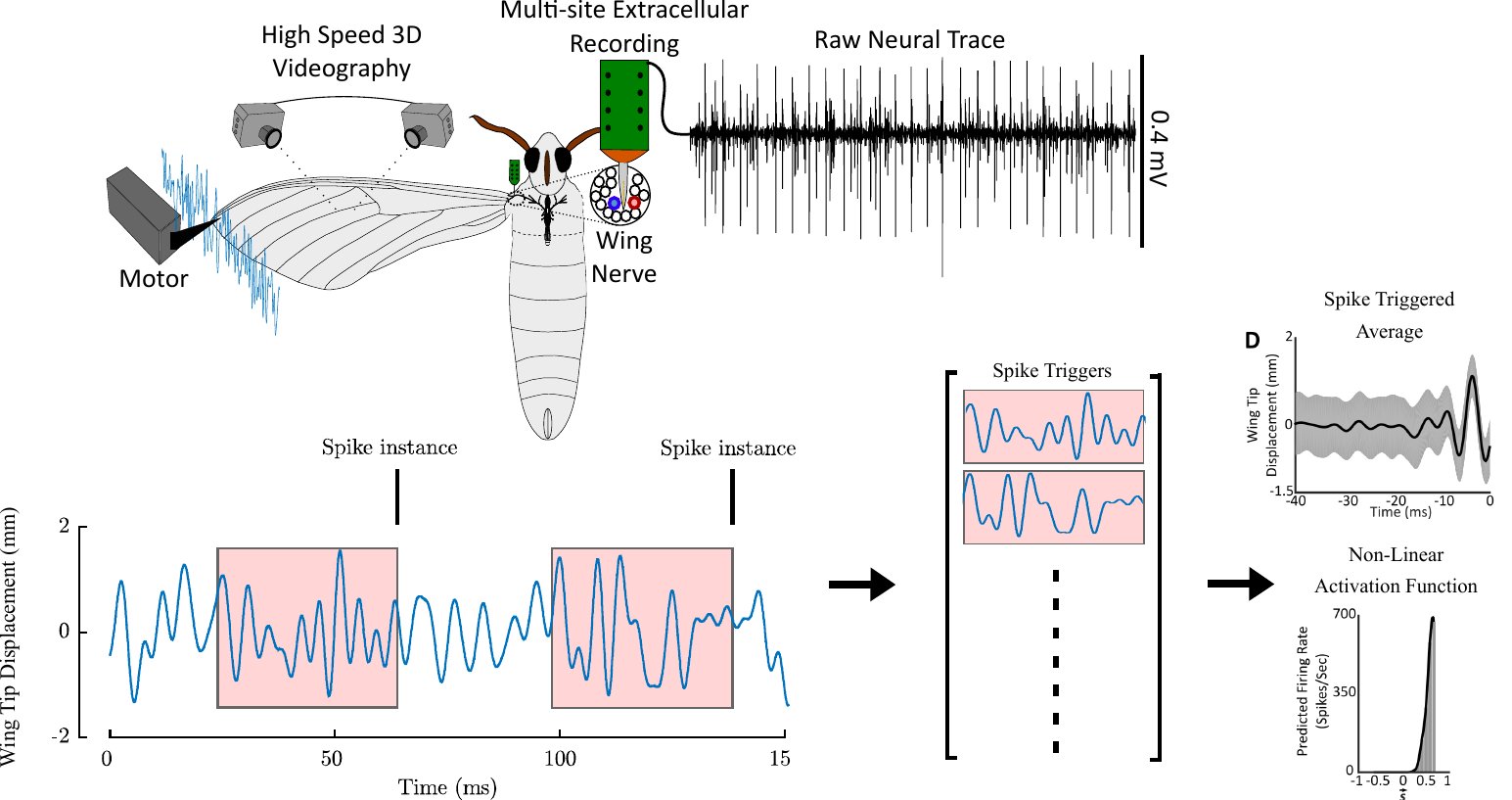}
\caption{
Pratt et al. \cite{pratt2017neural} characterized the neural response to strain of campaniform sensilla on hawkmoth wings by recording from the wing nerve and exciting the wing tip with band-limited white noise displacement. 
The strong feature selectivity of the spike triggered average and non-linear activation functions they observed were similar to observations in campaniform sensilla on the base of halteres of flies. 
}
\label{fig:Figure_Pratt}
\end{figure}

Campaniform sensilla on a wing fire action potentials in response to the time-history of mechanical forces they experience.
By recording from the wing nerve while stimulating the wing tip with a motor, one can characterize the stimulus that lead to neuronal firing.
A schematic diagram of the experimental setup is shown in Fig.~\ref{fig:Figure_Pratt}, and details of the experiment are found in \cite{pratt2017neural}. 

\subsubsection{Fitting STA \& NLA functions}
To summarize the responses of campaniform sensilla to mechanical stimulus, we compute the spike triggered average (STA) and a nonlinear activation (NLA) function to fit the experimental recordings.

The STA is approximated as a function of time $t$ before a spike at $t=0$,
\begin{align}
\text{STA}(t, f_\text{STA}, a,b) & =      \cos \big(  f_\text{STA}  (t+a)  \big)          \exp \left(  \frac{ - (t+a)^2}{b^2 }  \right), \label{eq:STA} \\
\end{align}
where $f_\text{STA}$ is the STA frequency, $a$ is the delay, and $b$ is the width.

\begin{figure}[h!]
 \centering\includegraphics[width=0.5\textwidth]{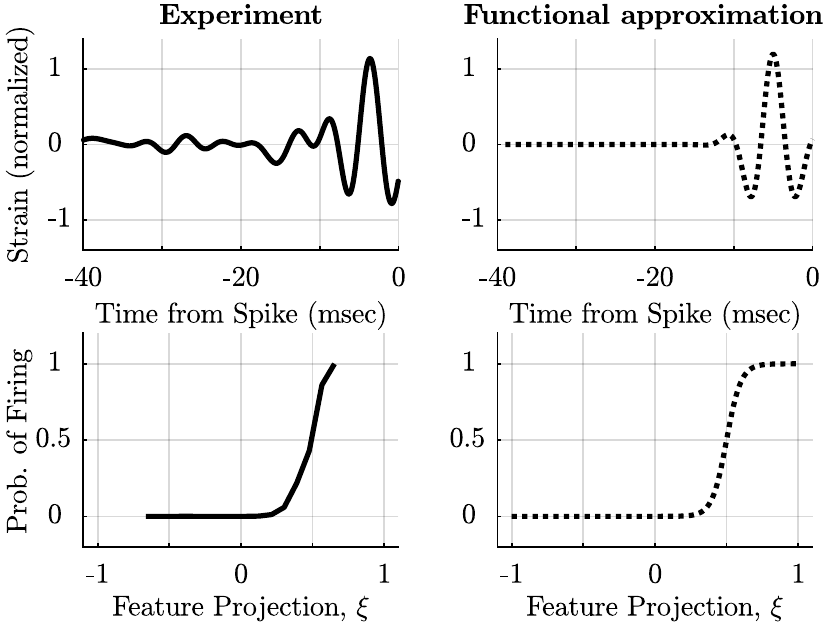}
\caption{
We approximated a typical experimentally observed feature and non linear activation functions from \cite{pratt2017neural} with parameterized functions according to equations \ref{eq:STA} and \ref{eq:NLD}. 
}
\label{fig:Figure_functionalApproximation}
\end{figure}

The strain $\epsilon(x, y, t)$ is convolved with the STA to obtain the strain projection on this STA feature $\xi(x, y, t)$. 
Next, $\xi(x, y, t)$ is mapped through a nonlinear activation function
\begin{align}
\text{NLA}(\xi,c,d) & = \frac{ 1}{  1 + \exp \left(  - c (\xi- d) \right) },  \label{eq:NLD}
\end{align}
where $c$ determines the slope and $d$ is the position of the function at half maximum.
Fig.~\ref{fig:Figure_functionalApproximation} shows the experimental STA and NLA as well as their best functional approximations.

\subsubsection{A probabilistic firing model} 
\begin{figure}[h!]
 \centering\includegraphics[width=0.75\textwidth]{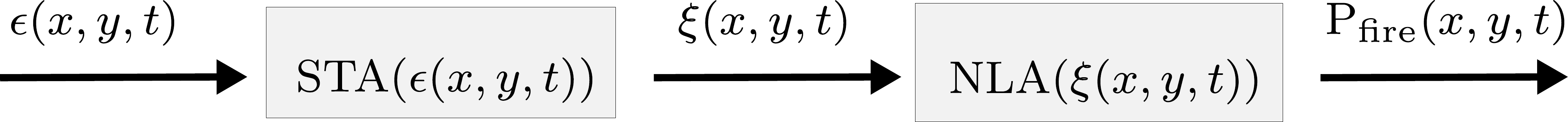}
\caption{ The strain $\epsilon(x, y, t)$ is converted to probability of firing by the two step process.
First we apply the discrete convolution of the strain with the STA and second, the NLA takes the resulting feature projection $\xi$, and outputs the probability of firing. 
}
\label{fig:Figure_STA_NLD_flowchart}
\end{figure}

We use the STA and NLA functions to transform the raw strain data from our structural simulation into probability of firing through a two-step process (Fig.~\ref{fig:Figure_STA_NLD_flowchart}). 
First, we apply a discrete convolution to the strain with the STA to obtain $\xi$,

\begin{align}
\xi(x,y,t)   & = \frac{1}{C_\xi} \sum_{ \tau = -39}^0   \epsilon(x,y,t-\tau)  \cdot \text{STA} (\tau).    \label{eq:STAapply}  
\end{align}

Second, we input $\xi$ into the NLA (\eqref{eq:NLAapply}).  
\begin{align}
P_\text{fire}(x,y,t)  & = \text{NLA} ( \xi(x,y,t) ).    \label{eq:NLAapply} 
\end{align}
The output is the probability of firing an action potential, which we define as the neural encoded strain.
Here $C_\xi$ is a constant to normalize the probability of firing and is determined by taking the maximum non-normalized $\xi$ over all sensors. 
$t$ is time in milliseconds. 
The probability of firing over the wing over time will form our data matrix $\mathbf{X}$ in the next section. \\

\nomenclature[RP]{$P_\text{fire} $}{Probability of neuron firing an action potential}%

%

\newpage
\subsection{Formulating the classification task} \label{ss:classification}

The strain data computed by the structural model as data to formulate a classification task. 
This section details how we define the training data, fit the classifiers, and assess the performance of classifiers on validation data.

\subsubsection{Constructing the data matrix} 

To construct build a classifier, we first define the training data.
The data matrix $\mathbf{X}$ comprises vectorized strain data from the two classes, flapping alone and flapping with rotation. 
Each row in $\mathbf{X}$ is data from a single sensor. 
For the full sensor set, $\mathbf{X}$ has 1326 rows.

\begin{align}
\mathbf{X}= \left[ \begin{array}{cccccc}
       \vdots 			 		&         	& \vdots 	       &   \vdots 			 		&         	& \vdots 	\\[0.3em]
      \mathbf{x^F_1} 	 			& \ldots 	& \mathbf{x^F_{k}} &  \mathbf{x^R_1} 	 			& \ldots 	& \mathbf{x^R_{k}}\\[0.3em]
      \vdots  &   \sunderb{6em}{ \text{Flapping}}  	&\vdots  &   \vdots   &  \sunderb{6em}{\text{With rotation} }  & \vdots 
     \end{array}  \right],
\end{align}

where $\mathbf{x^F}$ is strain data from flapping alone, $\mathbf{x^R}$ is strain data from flapping with rotation, and $k$ is the time index.

\subsubsection{Training and validation}
\begin{figure}[h]
    \centering
        \includegraphics{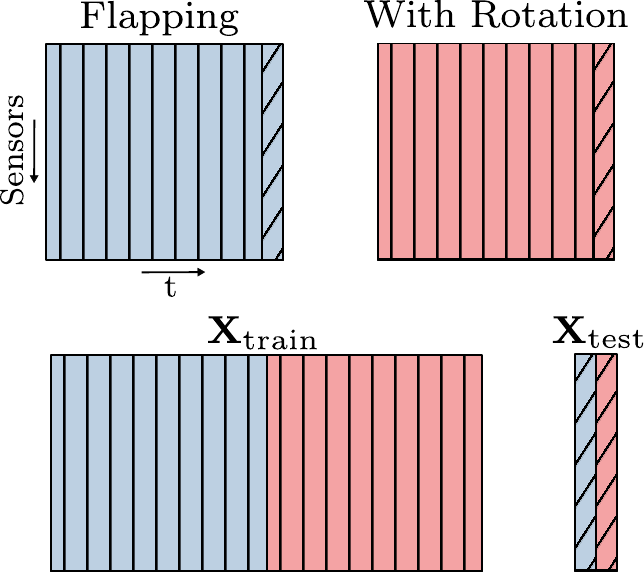}
        \caption{Construction of the training and validation data. The first 90\% of data snapshots of each class is assembled into $\mathbf{X_\text{train}}$, the last 10\% are assembled into $\mathbf{X_\text{test}}$. }
     \label{fig:Xtraintest}
\end{figure}

The training data is made from the first 90\% of snapshots for each class, which are assembled into $\mathbf{X_\text{train}}$.
The last 10\% of snapshots of data, taken from an epoch of the simulation after the training data, make $\mathbf{X_\text{test}}$.
We use $\mathbf{X_\text{train}}$ to fit the classifier, which will be used to predict the class for each of the snapshots in $\mathbf{X_\text{test}}$.  
This construction is shown schematically in Fig.~\ref{fig:Xtraintest}.

\subsubsection{Building linear classifiers} 

To classify the data, we use linear discriminant analysis (LDA) on the training data.
LDA computes a vector $\mathbf{w}$, where the data projected onto $\mathbf{w}$ is linearly separable.
In particular, we solve for
\begin{align}
\mathbf{w} = \argmax_{\mathbf{w'} }  \frac{   \mathbf{w'}^T \mathbf{S_B} \mathbf{w'}   }  {   \mathbf{w'}^T \mathbf{S_W} \mathbf{w'}   }  \label{eq:lda},
\end{align}
where $\mathbf{S_B}$ and $\mathbf{S_W}$ are the between-class and within-class scatter matrices, respectively. They are computed from the training data as follows:
\begin{align}
\mathbf{S_W}  &= \sum_{j=1}^c \sum_{ i \in c_j} (\mathbf{X}_{\text{train}, i} - \mu_j ) (\mathbf{X}_{\text{train}, i} - \mu_j )^T, \\
\mathbf{S_B} &= \sum_{j=1}^c  N_j (\mu_j - \mathbf{\mu}) (\mu_j - \mathbf{\mu} )^T. \\
\end{align}
Here $c$ is the total number of classes, $c_j$ are all the observations in the $j^{th}$ class, $\mu_j$ is the centroid of class $j$, and $\mathbf{\mu}$ is the centroid of all the training data.

To solve  \eqref{eq:lda}, we solve for the eigendecomposition of $\mathbf{S_W}^{-1} \mathbf{S_B}$, and $\mathbf{w}$ is the eigenvector corresponding to the largest eigenvalue:
\begin{align}
\mathbf{S_W}^{-1} \mathbf{S_B}  \mathbf{ w}   = \lambda \mathbf{ w }.
\end{align}

Note that the above linear classifier can be trained if the number of examples in the training set is at least $q + c$ to guarantee that $\mathbf{S_W}$ is not singular, where $q$ is the number of sensors and $c$ is the number of classes to be classified. 

Next, the data is projected onto $\mathbf{w}$ 
\begin{align}
\eta  = \mathbf{w}^T \mathbf{X}_\text{train},
\end{align}
and we apply a threshold in $\eta$ that separates the two classes.

\begin{figure}[h]
\centering
\includegraphics{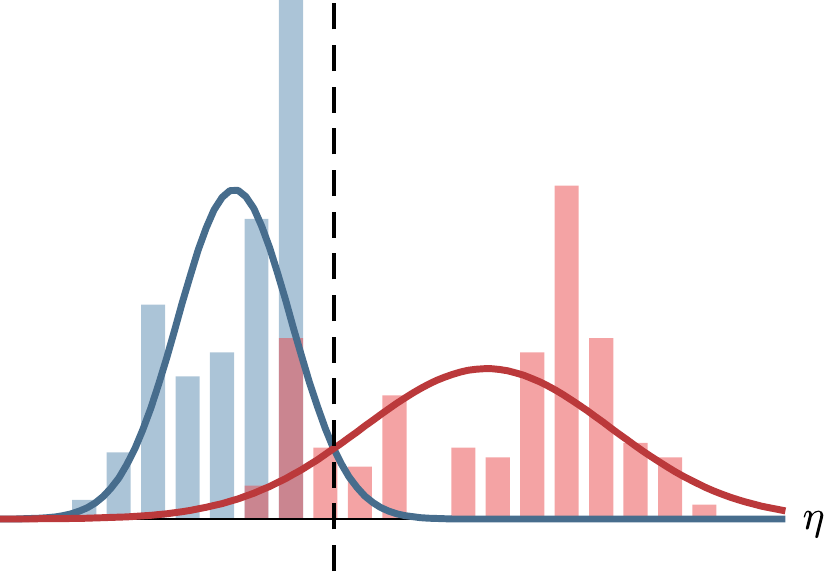}
\caption{The threshold between the two classes is determined by the intersection of the Gaussian probability density functions belonging to each class. }
\label{fig:threshold}
\end{figure}

For each class, we determine the mean $\mu_j$ and standard deviation $\sigma_j$. Assuming a gaussian distribution, we solve for where the two Gaussian probability density functions (\ref{eq:gaussian}) intersect in between the two means (Fig. \ref{fig:threshold}). If there is no intersection in between $\mu_1$ and $\mu_2$, we take the threshold to be the geometric middle $\frac{\mu_1 + \mu_2}{2}$. The Gaussian probability density function is given by: 

\begin{align}
G_j = \frac{1}{ \sqrt{2 \pi \sigma_j^2}}  \exp \left( - \frac{( x-\mu_j)^2}{2 \sigma_j^2} \right). \label{eq:gaussian}
\end{align}

\subsubsection{Classifier validation}

After training classifier by fitting $\mathbf{w}$, this classifier is applied to and assessed on the validation data $\mathbf{X}_{\text{test}}$,
\begin{align}
\eta_\text{test}  = \mathbf{w}^T \mathbf{X}_\text{test},
\end{align}
and the same threshold fit to the training data is applied to assign each test sample to a category. 
The validated accuracy is computed by comparing these categories to the known classes from the simulation.

\subsection{Sparse sensor placement optimization} \label{ss:sspoc}

This section describes our approach to learn a small handful of sparsely placed sensors that perform the body rotation classification.
We summarize the sparse sensor placement for optimal classification (SSPOC,~\cite{brunton2016sparse}).
Next, we describe two ways the sparse optimization we used in this paper have been extended from \cite{brunton2016sparse}.
First, we truncate selective singular vector features based on the discriminant vector $\mathbf{w}$.
And second, we use an elastic net penalty instead of a $\ell_1$ penalty.

\subsubsection{SSPOC} 

We take advantage of the observation that the high-dimensional data $\mathbf{x} \in \mathbb{R}^n$ may have a \emph{low-rank} representation:
\begin{align}
\mathbf{x} = \mathbf{\Psi}_r \mathbf{a},~~~~~~\mathbf{a} \in \mathbb{R}^r.
\end{align}
The goal of sparse sensor selection is to design a measurement matrix $\mathbf{C} \in \mathbb{R}^{q \times n}$ with a very small number of optimized measurements ($q \ll n$):
\begin{align}
\mathbf{y} = \mathbf{C} \mathbf{x} = \mathbf{C}  \mathbf{\Psi}_r \mathbf{a}.
\end{align} 
Further, we solve for $\mathbf{C}$ consisting of rows of the identity matrix, so that each sensor is a point measurement.

SSPOC~\cite{brunton2016sparse} is a sensor selection approach to find a solution $\mathbf{C}$ so that the linear discrimination between classes is achievable with the sparse measurements $\mathbf{y}$.
In briefly, we first reduce the dimensionality of $\mathbf{X}$ using the singular value decomposition (SVD):
\begin{align}
\mathbf{X} = \mathbf{\Psi} \mathbf{\Sigma} \mathbf{V}^T \approx  \mathbf{\Psi_r} \mathbf{\Sigma_r} \mathbf{V_r}^T,
\end{align}
where we take advantage of order of the singular values to truncate $\mathbf{\Psi}$, $\mathbf{\Sigma}$, and $\mathbf{V}$ to their first $r$ features.
The data in $\mathbf{X}$ may be projected to $\mathbb{R}^r$ using $\mathbf{\Psi_r}$,
\begin{align}
\mathbf{a} = \mathbf{\Psi_r}^T \mathbf{X}.
\end{align}
Next, we use LDA to solve for the discriminant vector $\mathbf{w}$ using $\mathbf{a}$ as the training data, so that the discriminant threshold is applied in $\eta$:
\begin{align}
\eta = \mathbf{w}^T  \mathbf{a} = \mathbf{w}^T \mathbf{\Psi_r}^T \mathbf{X}.
\end{align}

Finally, we solve for the sparse vector $\mathbf{s} \in \mathbb{R}^n$:
\begin{align}
\mathbf{s} = \argmin_{\mathbf{s'}} \left\lVert \mathbf{s'} \right\rVert_1, ~~~~~\text{subject to}~~ \mathbf{\Psi_r}^T \mathbf{s'} = \mathbf{w}, \label{eq:sspoc}
\end{align}
where $\mathbf{s}$ comprises mostly zeros, and the non-zero entries of $\mathbf{s}$ correspond to sensor locations and rows of the identity matrix selected for the measurement matrix $\mathbf{C}$.

Generally speaking, the number of sensors $q$ selected by this approach is approximately $r$, so the choice of $r$ determines the number of sensors desired.

\subsubsection{Singular value feature selection for SSPOC}

In the previous section, we described truncating the SVD basis $\mathbf{\Psi}$ to its first $r$ columns, corresponding to the $r$ largest singular values.
However, these first $r$ features may not necessarily be the ones supporting the largest separation between classes.
Here, we use an alternative criterion to select which columns of $\mathbf{\Psi}$ are used in the sparse optimization by re-weighing each according to the LDA discriminant vector.

Specifically, the singular values are re-weighted according to the magnitude of $\mathbf{w}$, and the $\rho$ largest entries of
$\mathbf{\Sigma_r} |\mathbf{w}|$
determine the column of $\mathbf{\Psi}$ that form a new truncated basis $\mathbf{\Psi}_\rho$.
It follows that we solve \eqref{eq:sspoc} using $\mathbf{\Psi}_\rho$ and $\mathbf{w}_\rho$, which produces approximately $\rho$ sensors.

\subsubsection{Sparse optimization with an elastic net penalty} 

We observed that solutions to \eqref{eq:sspoc} using convex optimization tools sometimes do not converge to optimal solution.
Therefore, in this paper we use a related optimization using an elastic net penalty, which balances the ratio of penalty for the $\ell_1$ and $\ell_2$ norms of $\mathbf{s}$:
\begin{align}
\mathbf{s} &=  \argmin_{\mathbf{s'}} \  \alpha ||\mathbf{s'}||_1  + (1-\alpha) ||\mathbf{s'}||_2 , \ \mathrm{subject} \  \mathrm{to} \ \mathbf{\Psi_\rho^T s'} = \mathbf{w_\rho}.\label{eq:elasticnet}
\end{align}
For this paper, we use $\alpha = 0.9$.

\nomenclature[GP]{$\mathbf{  \Psi}$}{A basis in which the data $\mathbf{X}$ may be represented, for instance the singular value decomposition (SVD) basis}%
\nomenclature[GPr]{$\mathbf{  \Psi_r }$}{Matrix of left singular vectors truncated to the first $r$ columns}%
\nomenclature[GPt]{$\mathbf{  \Psi_\rho }$}{Matrix of left singular vectors used for sensor selection for classification}%
\nomenclature[GSr]{$ \mathbf{ \Sigma}$}{Matrix of singular values}%

\nomenclature[Rw]{$ \mathbf{ w } $}{ Linear discriminant vector }%
\nomenclature[Ga]{$ \eta$}{Decision variable}%
\nomenclature[Rs]{$\mathbf{s} $}{Vector of mostly zeros the same size as $\mathbf{x}$, where the nonzero element correspond to desired sparse sensor locations}%
\nomenclature[Rq]{$q$}{Number of sensors}%

\subsection{Sigmoidal fit to classification accuracy}  \label{ss:sigmoid}
We observed that the validated accuracy $A$ depends on the number of sensors $q$ in a sigmoidal relationship, so we fit the results shown in Figs.~\ref{fig:Figure_R1} and \ref{fig:Figure_R2} with a sigmoidal function with 3 parameters:
\begin{align}
A(q) = \frac{  \frac{1}{2} + c_1 }{  1 + \exp(  - \frac{ q - c_2}{c_3} )  }.  \label{eq:sigmoidFit} 
\end{align}

To summarize these curves over different neural encoders, the results in Fig.~\ref{fig:Figure_R3} present the fewest number of sensors required for classification at 75\% accuracy. 
This number is determined by solving for $q$ at which $A(q)$ crosses $0.75$.

 \begin{figure}[h]
\centering
\includegraphics{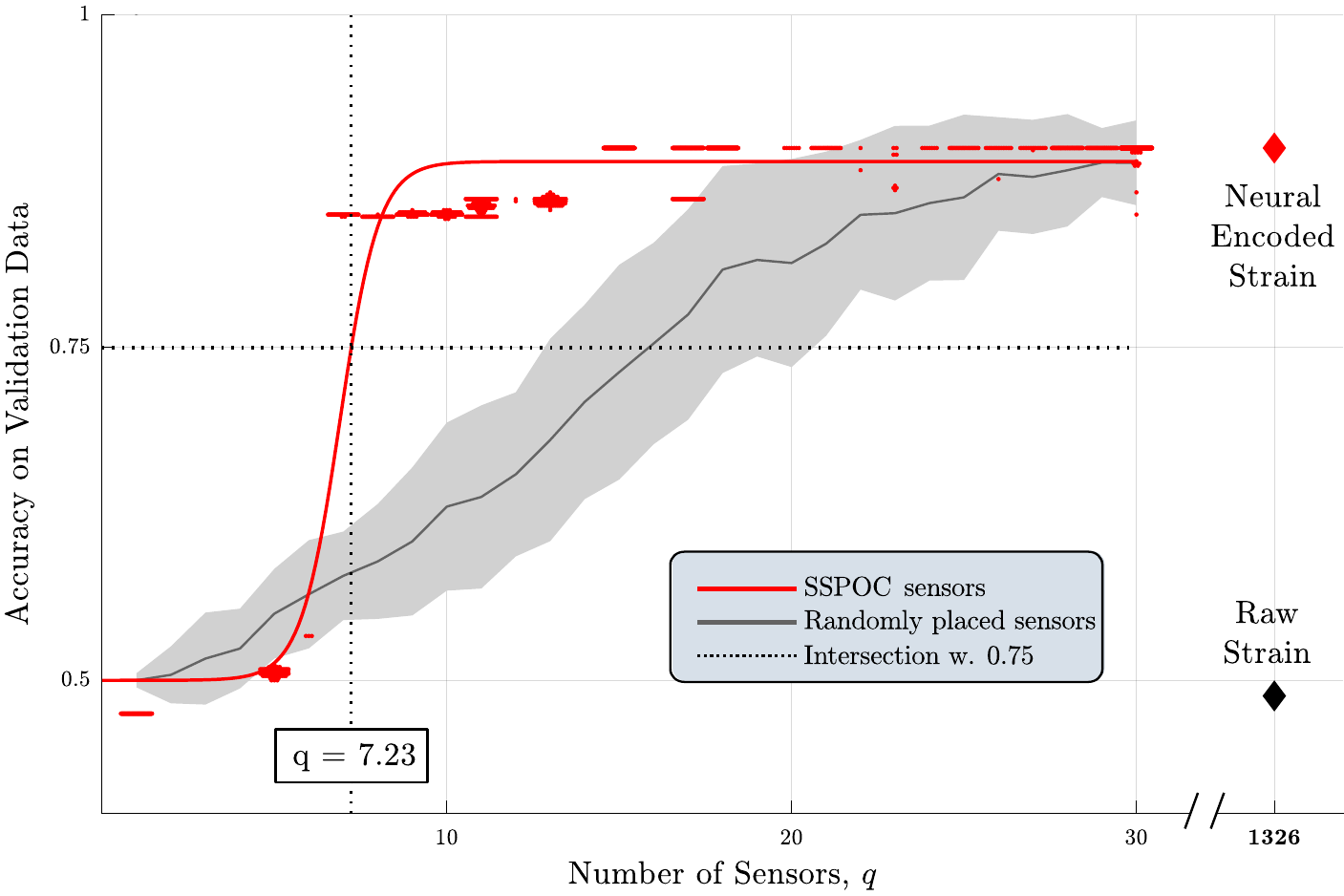}
\caption{
We determined the number of sensors required for good classification by fitting a sigmoid (eq. \ref{eq:sigmoidFit}) to the classification accuracy versus the number of sensors. 
The dashed line shows the sigmoidal fit with constants $[c_1,c_2,c_3] = [0.378, 6.904, 0.583] $, intersecting with the 0.75 accuracy line at q = 7.29. 
  }
\label{fig:Figure_S6}
\end{figure}

\newpage
\subsection{Supplementary results figures} \label{ss:supplresults}

\end{document}